\documentclass[11pt,a4paper,reqno]{amsart}%
\usepackage{amsthm,amsmath,amsfonts,amssymb,xcolor,amsxtra,appendix,bookmark,dsfont,bm,mathrsfs}
\usepackage[margin=38mm]{geometry}
\usepackage[latin1]{inputenc}
\usepackage{color} 
\usepackage{amssymb}\usepackage{graphicx}
\usepackage{tikz}
\usepackage{hyperref}
\theoremstyle{plain}
\theoremstyle{remark}

\allowdisplaybreaks

\title[Wave turbulence and collective behavior dynamics]{Wave turbulence and collective behavior models for  wave equations with short- and long-range interactions}

\author[A. Aceves]{Alejandro Aceves
}
\address{Department of Mathematics, Southern Methodist University, Dallas, Texas 75275, USA}
\email{aaceves@mail.smu.edu} 
\thanks{A.~A. is  funded in part by  NSF RTG Grant DMS-1840260 and NSF-DMS-1909559}

\author[R. Alonso]{Ricardo Alonso
}
\address{Department of Mathematics, Texas A \& M University at Qatar, PO Box 23874, Education City Doha, Qatar}
\email{ ricardo.alonso@qatar.tamu.edu} 
\thanks{R.~A. is funded in part by Bolsa de Produtividade em Pesquisa CNPq (303325/2019-4).}

\author[M.-B. Tran]{Minh-Binh Tran}
\address{Department of Mathematics, Southern Methodist University, Dallas, Texas 75275, USA}
\email{minhbinht@mail.smu.edu} 
\thanks{M.-B. T. is  funded in part by  the  NSF Grant DMS-1854453, NSF RTG Grant DMS-1840260,  URC Grant 2020, Humboldt Fellowship,  Dedman College Linking Fellowship, NSF CAREER  DMS-2044626.}

\begin{document}
\date{\today}

\begin{abstract} In this work, we discuss a situation which could lead to both wave turbulence and collective behavior kinetic equations. The wave turbulence kinetic models appear in the kinetic limit when the wave equations have local differential operators. Viewing wave equations on the lattice as chains of anharmonic oscillators and replacing the local differential operators (short-range interactions) by non-local ones (long-range interactions), we arrive at a new Vlasov-type kinetic model in the mean field limit under the molecular chaos assumption reminiscent of models for collective behavior in which anharmonic oscillators replace individual particles.
\end{abstract}

\maketitle
\centerline{\it Dedicated to the  85th birthday of Professor Roland Glowinski}
 \tableofcontents

\medskip

\section{Introduction}\label{Sec:Main} 
 Having the origin in the works of  Peierls \cite{Peierls:1993:BRK,Peierls:1960:QTS},  Hasselmann \cite{hasselmann1962non,hasselmann1974spectral},  Benney-Saffman-Newell \cite{benney1969random,benney1966nonlinear},  Zakharov \cite{zakharov2012kolmogorov},   wave turbulence (WT) theory describes the dynamics of weakly nonlinear and dispersive waves (classical or non-classical)  out of thermal equilibrium. Even though  wave fields describing the processes of random wave interactions  in nature   are  enormously diverse, a common mathematical framework can be used to model the dynamics of spectral energy transfer  in both quantum or classical wave systems. In this mathematical framework, the probability density functions associated with weakly nonlinear wave interactions are solutions of wave kinetic (WK) equations. Over the years,  WK equations have been shown to play  important roles in a vast range of physical applications, as discussed in the books  \cite{zakharov2012kolmogorov,Nazarenko:2011:WT}.  {We also mention closely related kinetic models developed when the interest focuses in the interaction of particles and oscillators, see for example \cite{Parris}. }
 
In addition, since the  realization of Bose-Einstein condensation (BEC) in trapped atomic vapors of $^{23}$Na \cite{davis1995bose}, $^{87}$Rb \cite{WiemanCornell} and $^7$Li \cite{bradley1995evidence}, a period of intense theoretical and experimental   research  has been initiated.  A theoretical Quantum Kinetic (QK) theory, which takes into account the coupled non-equilibrium dynamics of both  the thermal cloud of the Bose gas and the BEC  under investigation, is needed to support the experimental results.  Although being used to described different physical phenomena, QK kinetic equations are quite similar with WT ones \cite{PomeauBinh,tran2020boltzmann,tran2021thermal}. During the last few years, there has been a growing interests in rigorously understanding those kinetic equations. 
Starting with the pioneering work  of Lukkarinen and Spohn \cite{LukkarinenSpohn:WNS:2011}, there have been a lot of recent works in 
in rigorously deriving WK equations (see,  for instance \cite{ampatzoglou2021derivation,buckmaster2019onthe,buckmaster2019onset,collot2019derivation,collot2020derivation, deng2019derivation,deng2021full,dymov2019formal,dymov2019formal2,dymov2020zakharov,dymov2021large,staffilani2021wave} and the references therein). The analysis of WK and QK equations is also a topic of current interest. We list here an incomplete list and refer to the references therein for more detailed descriptions of the literature   \cite{AlonsoGambaBinh,AlonsoGoudon,ArkerydNouri:2012:BCI,ArkerydNouri:AMP:2013,ArkerydNouri:2015:BCI,cai2018spatially,Goudon,EscobedoMischlerValle:HBI:2003,Binh9,EscobedoVelazquez:2015:FTB,EscobedoVelazquez:2015:OTT,GambaSmithBinh,germain2020optimal,JinBinh,Lu:2004:OID,Lu:2005:TBE,Lu:2013:TBE,ToanBinh,nguyen2017quantum,Binh1,SofferBinh2,soffer2019energy,CraciunSmithBoldyrevBinh}.

Collective behavior of self-propelled particles such as swarming of bacteria,  schooling of fishes, flocking of birds and mobile agents, appears in many contexts \cite{blondel2009krause,couzin2005effective,cucker2007mathematics,cucker2007emergent,during2009boltzmann,hegselmann2002opinion,motsch2011new,Orsonga,vicsek1995novel,zavlanos2011graph}.
How and when  clusters emerge, and what type of rules of engagement in influences clusters are among the questions that have been attracted the attention of scientists for decades. Over the last years, there have been growing interests in the mathematical community in studying those models rigorously. We list here only a few of those works and refer the readers to the references therein as the list is quite incomplete \cite{albi2019vehicular,bellomo2012mathematical,Bolley,caponigro2015sparse,carrillo2014derivation,degond2015phase,duan2010kinetic,ha2009emergence,ha2008particle,he2021game,karper2015hydrodynamic,ko2020asymptotic,morales2019flocking,motsch2014heterophilious,shvydkoy2017eulerian,tadmor2021mathematics,tadmor2014critical}.

In this paper, we discuss a  connection between wave turbulence and collective behavior kinetic models. Starting from the weakly nonlinear wave equation on the lattice, it has been showed \cite{LukkarinenSpohn:WNS:2011,staffilani2021wave} that the wave kinetic equations  can be derived rigorously, under suitable assumptions on the randomization of the initial condition and the wave equations. This procedure is summarized in Section \ref{WaveT} for a wave equation with a quadratic nonlinearity and the kinetic equation under consideration is the 3-wave kinetic equation  \eqref{WaveTurbulence}. However, by using other type of nonlinearities, we could arrive at different wave kinetic equations, with the same procedure. It is well-known in the physical community that chains of anharmonic oscillators, such as the Fermi-Pasta-Ulam-Tsingou (FPTU),  also exhibit  collective behaviors, if the interactions are long-range  (see for instance \cite{christodoulidi2016dynamics,christodoulidi2014fermi}). As wave equations on the lattice could also be viewed as chains of anharmonic oscillators, by replacing the local differential operators (for example, the Laplace or Biharmonic operators) in the wave equations by operators that describe long-range interactions (for example, the fractional Laplace operators), we could expect to obtain  models that exhibit  collective behaviors.  One of the key difference between our chains of anharmonic oscillators and models for consensus, flocking and swarming \cite{cucker2007mathematics,cucker2007emergent,ha2008particle} is that, in our case, the system under consideration will need to ``label" the location of the oscillator in the lattice.  This is done by adding a new kinetic variable in the density distribution of anharmonic oscillators representing such location.

In Section \ref{Collective}, we derive formally three Vlasov-type kinetic equations \eqref{Vlasov1}, starting from a wave equation whose differential operator is a fractional Laplacian.  Here, the density  $g$ is not only a function of the position $r$ (of the oscillation), the velocity $v$ (of the oscillation), the time $t$ variables but also of an additional continuous variable $x$, which ``labels'' the location of the anharmonic oscillator in the lattice.  This framework is reminiscent of polyatomic and multicomponent models that add a kinetic variable to differentiate species, see for example \cite{Bisi} and references therein.

%Due to the fact that the particles cannot be considered as indistinguishable, in \eqref{Vlasov1}-\eqref{Vlasov3}, 
% In comparison with the Vlasov equation obtained from the C-S model (see \cite{ha2008particle}), the equations \eqref{Vlasov1}-\eqref{Vlasov3} have an extra variable $x$, which could be seen as a ``label'' for the different particles, besides   the standard space $r$,  velociy $v$,  time $t$ variables.

Let us mention, however,  that the concept of ``label'' has been previously introduced in consensus models in the work of Biccari, Ko and Zuazua \cite{biccari2019dynamics}, which considers ``networked consensus models'' and has an inspiration from the previous work of Kawamura \cite{kawamura2014kuramoto} on the nonlocal Kuramoto-Sakaguchi equation, where the ``label'' is indeed as a location vector. A different form of the  Kuramoto-Sakaguchi equation is used in \cite{amadori2017global}, where the distribution function $f$ depends on the state  $\theta$ and the natural frequency $\Omega$, which represents the ``label''.

{These type of Vlasov-type kinetic equations are recent and increasingly important in the literature and their mathematical properties are open for investigation.}  Due to the scope of our paper, in Section \ref{Collective}, we only focus on an example of a linear wave equation.

In considering collective behavior, we are reminded of the surprising recurrence result in FPUT, that countered the expected thermalization. It is perhaps the nearest neighbor coupling (i.e. local interaction as in classical random walk) combined with nonlinearity that triggers this behavior. Alternatively, in the canonical model for global coupling, the Kuramoto model \cite{kuramoto1975international}, at sufficiently high coupling strength, collective synchronous behavior emerges, which overcomes the expected deviations from the natural frequency of the ideally identical oscillators. The question is, if these observations are representative in our newly obtained Vlasov-type kinetic equations. A perhaps subtle but important difference between the Kuramoto model and wave-like and corresponding Sch\"odinger-like models, is that in the first one the state variable referred to as an angle is real, whereas in particular in applications on electromagnetism and quantum mechanics, in the second case the state variable is complex. In this second case recent research in quantum mechanics and photonics suggests the importance of long range (global) coupling of for example photonic resonators of fiber amplifiers, to enhance coherence. In fact, it is coherence that one views as the order parameter than in the Kuramoto models measures the degree of synchronization.

{\bf Acknowledgements.}  We dedicate this paper to the Special Issue of the Journal in honor of Professor Roland Glowinski. 
M.-B. T would like to express his attitude toward Professor Roland Glowinski for his constant guidance, support and friendship over the years. M.-B.T would also like to thank   Dongnam Ko for the fruitful remarks on Collective Behavior Theory and the explanations of the work \cite{amadori2017global,biccari2019dynamics,kawamura2014kuramoto}.
\section{Wave turbulence kinetic models for discrete nonlinear wave equations with short-range interactions}\label{WaveT}

Let us first start with a nonlinear wave equation, with a quadratic nonlinearity. However, our discussion could be extended to cubic and higher order nonlinearities. 
\begin{equation}
	\label{QuadraticNLS}\begin{aligned}
		&	\frac{\partial^2\psi}{\partial t^2}(x,t) \ + \ \mathcal{L} \psi(x,t)   \ + \lambda\psi^2(x,t) \ = \ 0,\\
		\psi(x,0) \ & = \ \psi_0(x), \ \ \frac{\partial\psi}{\partial t}(x,0) \ = \ \psi_1(x),\end{aligned}
\end{equation}
for $x$ being on the torus $[0,1]^d$, $t\in\mathbb{R}_+$, %$\omega_0\in\mathb%b{R}$ 
%is some real %constant,
$\lambda$ is a small constant describing the smallness of the nonlinearity. We suppose that the interactions are short-range, which is, the operator $\mathcal L$ is a standard local differential operator, for instance, when $\mathcal{L}=\Delta$, we obtain the Klein-Gordon equation and when $\mathcal{L}=-\Delta^2$, we obtain the beam wave equation. Similar with \cite{staffilani2021wave}, we introduce the finite volume mesh, namely
\begin{equation}
	\label{Lattice}
	\Lambda \ = \ \Lambda(D) \ = \ \left\{0,\frac{1}{2D+1} \dots,\frac{2D}{2D+1}
	\right\}^d, 
\end{equation}
for some constant $D\in\mathbb{N}$. As we will work on the Fourier transform, we define the mesh size of the frequency space to be \begin{equation}\label{Mesh} h=\frac{1}{2D+1}.\end{equation} 
We  follow \cite{staffilani2021wave} and introduce discretized equation 
\begin{equation}
	\label{LatticeDynamics}\begin{aligned}
		\partial_{tt}\psi(x,t) \ = &  \ \ -\sum_{y\in  \Lambda} O_1(x-y)\psi(y,t) \ - \ \lambda ({\psi(x,t)})^2, \\
		\psi(x,0) \ & = \ \psi_0(x), \ \partial_{t}\psi(x,0) \  = \ \psi_1(x), \ \forall (x,t)\in\Lambda\times \mathbb{R}_+,\end{aligned}
\end{equation}
in which $O_1(x-y)$ is the finite difference operator obtained from the continum operator $\mathcal L$. We now introduce  the discrete Fourier transform 
\begin{equation}
	\label{Def:Fourier}\hat  \psi(k)=h^d \sum_{x\in\Lambda} \psi(x) e^{-2\pi {\bf i} k\cdot x}, \quad k\in \Lambda^* = \Lambda^*(D) =  \{-D,\cdots,0,\cdots,D\}^d.
\end{equation}
At the end of this standard procedure, \eqref{LatticeDynamics}  can  be rewritten in the Fourier space as a system of ODEs 
\begin{equation}
	\label{LatticeDynamicsFourier}\begin{aligned}
		\partial_{tt}\hat\psi(k,t) \ = \  &  - \big(\bar\omega(k)\big)^2\hat\psi(k,t)\ 
	- \lambda \sum_{k=k_1+k_2;k_1,k_2\in\Lambda^*}{\hat\psi(k_1,t)}\hat\psi(k_2,t),\\
		\hat\psi(k,0) \ & = \ \hat\psi_0(k) ,\ \ \ \partial_{t}\hat\psi(k,0) \ = \ \hat\psi_1(k). \end{aligned}
\end{equation}
In the beam wave case, $\mathcal{L}=-\Delta^2$ and the  dispersion relation  takes the  discretized form (see \cite{rumpf2021wave})
\begin{equation}
	\label{NearestNeighbordA:Bretherton}
	\bar\omega(k) \ = \ \sin^2(2\pi h k^1) + \cdots + \sin^2(2\pi h k^d),  
\end{equation}
with $k=(k^1,\cdots,k^d)$. Later, we will also need the rescaled dispersion relation
\begin{equation}
	\label{omegaB}
	\omega(k) \ = \ \sin^2(2\pi k^1) + \cdots + \sin^2(2\pi  k^d).  
\end{equation}

% We also randomize the initial condition as follows
%
%\begin{equation}
%	\begin{aligned}\label{InitialCondB}
%		\hat\psi_0(k,\varepsilon,\tilde\varepsilon) \ = \ &\frac{1}{\omega^b(k)} \sum_{\xi\in\Lambda_B}\hat\psi_0^2(k-\xi/\epsilon,{\varepsilon})\hat\psi_0^1(\xi,\tilde{\varepsilon})\\
%			\hat\psi_1(k,\varepsilon,\tilde\varepsilon) \ = \ &\frac{\omega^b(k)}{{\bf i}} \sum_{\xi\in\Lambda_B}\hat\psi_0^3(k-\xi/\epsilon,{\varepsilon})\hat\psi_0^1(\xi,\tilde{\varepsilon}),
%	\end{aligned}
%\end{equation}
%where the random variables $\hat\psi_0^2(k-\xi/\epsilon,{\varepsilon}),\hat\psi_0^3(k-\xi/\epsilon,{\varepsilon})$ are i.i.d and follow the probability law of $W^b_k(t,\varepsilon)$. The random variables $\hat\psi_0^1(\xi,\tilde{\varepsilon})$ defined on $(\Omega',{\bf{P}}')$ are i.i.d and  independent of $W^b_k(t,\varepsilon)$.

We define the inverse Fourier transform 
\begin{equation}\label{Def:FourierInverse}
	f(x)= \sum_{k\in\Lambda_*} \hat  f(k) e^{2\pi {\bf i} k\cdot x},
\end{equation}
as well as the shorthand notations
\begin{equation}
	\label{Shorthand1}
	\int_{\Lambda}\mathrm{d}x \ = \  h^d\sum_{x\in\Lambda},\ \ \ \ 
	\langle f, g\rangle \ = \ h^d \sum_{x\in\Lambda}f(x)^* g(x),\ \ \ \ \langle x\rangle \ = \ \sqrt{1+|x|^2}, \ \ \forall x\in\mathbb{R}^d,
\end{equation}
where   if $z\in \mathbb{C}$, then ${z}^*$ is the complex conjugate. We also denote
\begin{equation}
	\sum_{k\in\Lambda^*} \ = \ \int_{\Lambda^*}\mathrm{d}k.
\end{equation}
In addition, for any $N\in\mathbb{N}\backslash\{0\}$, similar with \cite{staffilani2021wave}, we define the delta function $\delta_N$ on $(\mathbb{Z}/N)^d$ as
\begin{equation}
	\label{Def:Delta}\delta_N(k) = |N|^d\mathbf{1}(k \mbox{ mod } 1 \ = \ 0), \ \ \ \forall k\in (\mathbb{Z}/N)^d,
\end{equation}
in which the sub-index $N$ is commonly omitted and written as
\begin{equation}
	\label{Def:Delta1}\delta(k) = |N|^d\mathbf{1}(k \mbox{ mod } 1 \ = \ 0), \ \ \ \forall k\in (\mathbb{Z}/N)^d.
\end{equation} 
Equation \eqref{LatticeDynamicsFourier} can now be expressed as a coupling system
\begin{equation}
	\label{LatticeDynamicsSystem}\begin{aligned}
		\partial _t {q}(k,t) \ & = \ {p}(k,t),\\
		\partial _t{p}(k,t) \ & = \ -\ \big(\bar\omega(k)\big)^2{q}(k,t)\\
		&\ \ \  \ -  \lambda\int_{(\Lambda^*)^2}\mathrm{d}k_1\mathrm{d}k_2\delta(k-k_1-k_2){q}(k_1,T){q}(k_2,t), \\
		q(k,0) \ & = \ \hat\psi_0(k), \ \ p(k,0) \  = \ \hat\psi_1(k), \ \ \ \ \forall (k,t)\in\Lambda^*\times \mathbb{R}_+,\end{aligned}
\end{equation}
which, under the transformation (cf. \cite{zakharov1967weak})
\begin{equation}\label{Tranformation1}
	a(k,t) \ =  \ {\bar\omega(k)}{q}(k,t)\ +  \ \frac{{\bf i}}{{\bar\omega(k)}}{p}(k,t),
\end{equation}
with the inverse
\begin{equation}\label{Tranformation2}\begin{aligned}
		{q}(k,T) \ & =  \ \frac{1}{2{\bar\omega(k)}}\Big[a(k)\ + 
		\ a^*(-k)\Big],\\
		{p}(k,T) \ & =  \ {\bf i}\frac{{\bar\omega(k)}}{2}\Big[-a(k)\ + 
		\ a^*(-k)\Big],
	\end{aligned}
\end{equation}
leads to the following system of ordinary differential equations
\begin{equation}
	\begin{aligned}\label{ODEs}
		\partial _t a(k,t)  \  =  \ & -{\bf i}\bar\omega(k)a(k,t) \ - \ {\bf i}\lambda\int_{(\Lambda^*)^2}\mathrm{d}k_1\mathrm{d}k_2\delta(k-k_1-k_2)\times\\
		&\  \times [8\bar\omega(k)\bar\omega(k_1)\bar\omega(k_2)]^{-1}\Big[a(k_1,t) \ + \ a^*(-k_1,t)\Big]\Big[a(k_2,t) \ + \ a^*(-k_2,t)\Big],\\
		a(k,0) \  = &\ a_0(k) \ = \ \frac12{\Big[{\bar\omega(k)}{q}(k,0)\ +  \ \frac{{\bf i}}{{\bar\omega(k)}}{p}(k,0)\Big]},  \forall (k,t)\in\Lambda^*\times \mathbb{R}_+.
	\end{aligned}
\end{equation}
%We also have the initial condition \begin{equation}
%	\begin{aligned}\label{InitialCondB:1}
%		a(k,0) \ = \ & \sum_{\xi\in\Lambda_B}[\hat\psi_0^2(k-\xi/\epsilon,{\varpi})+\hat\psi_0^3(k-\xi/\epsilon,{\varpi})]\hat\psi_0^1(\xi,\tilde{\varpi})
%		\\
%		\ =: \  & \sum_{\xi\in\Lambda_B}a'_0(k-\xi/\epsilon,{\varpi})a''_0(\xi,\tilde{\varpi}),
%	\end{aligned}
%\end{equation}
%with $a'_0(k-\xi/\epsilon,{\varpi}=\hat\psi_0^2(k-\xi/\epsilon,{\varpi})+\hat\psi_0^3(k-\xi/\epsilon,{\varpi})$ and $a''_0(\xi,\tilde{\varpi})=\hat\psi_0^1(\xi,\tilde{\varpi}).$

%As we are interested in only real solutions of \eqref{LatticeDynamics}, we deduce  $\hat{\psi}^*(k)=\hat{\psi}(-k)$, yielding  $a(-k)=a^*(k)$, then $\frac{{\bf i}\sqrt{2c_{r}}}{2}\Big[a(k,t) \ + \ a^*(-k,t)\Big]={{\bf i}\sqrt{2c_{r}}}a(k,t).$  
%
%
Let ${a},{a}^*$ denote the vectors $(a_k)_{k\in\Lambda^*}$, $(a_k^*)_{k\in\Lambda^*}$, and let us  set 
\begin{equation}\label{Hamiltonian}
	{H}({a},{a}^*) \ = \ {H}_1({a},{a}^*)  \ + \ \lambda{H}_2({a},{a}^*),
\end{equation}
with 
$${H}_1({a},{a}^*)  \ = \ \sum_{k\in\Lambda_B^*}\frac{1}{2}\bar\omega(k)|{a}_k|^2,$$
$${H}_2({a},{a}^*)  \ = \ \sum_{k,k_1,k_2\in\Lambda^*}\mathcal{W}(k,k_1,k_2)\delta(k-k_1-k_2)\Big[a(k_1,t) \ + \ a^*(-k_1,t)\Big]$$
$$\times\Big[a(k_2,t) \ + \ a^*(-k_2,t)\Big]{a}_{k}^*,$$ 
\begin{equation}
	\label{Kernel}
	\mathcal{W}(k,k_1,k_2) \ = \ [\bar\omega(k)\bar\omega(k_1)\bar\omega(k_2)]^{-1}, \ \ \  \mathcal{M}(k,k_1,k_2) \ = \ [8\bar\omega(k)\bar\omega(k_1)\bar\omega(k_2)]^{-1}.
\end{equation}

We then obtain the system
\begin{equation}\label{HamiltonODEs}
	\partial _t{a}_k \ = \ {\bf i}\frac{\partial{H}(a,a^*)}{\partial {a}_k^*}.\end{equation}
By defining 
\begin{equation}
	\label{ShortenNotationofA}
	\hat{a}(k,1,t) \ = \ a_k(t), \mbox{ and } \hat{a}(k,-1,t) \ = \ a^*_k(t),
\end{equation}
we rewrite the system \eqref{ODEs} as
\begin{equation}
	\begin{aligned}\label{StartPoint0}
		\partial _t \hat{a}(k,\sigma,t)\ & =  \ -{\bf i}\sigma\bar\omega(k)\hat{a}(k,\sigma,t) \mathrm{d}t\ - \ {\bf i}\sigma\lambda\sum_{\sigma_1,\sigma_2\in\{\pm1\}}\sum_{k_1,k_2\in\Lambda^*}\delta( \sigma k-\sigma_1k_1-\sigma_2k_2)\\
		& \ \   \times \mathcal{M}(k,k_1,k_2)\hat{a}(k_1,\sigma_1,t) \hat{a}(k_2,\sigma_2,t),\\
		\\
		\hat{a}(k,1,0) \ & = \ a_0(k), \ \ \ \ \forall (k,t)\in\Lambda^*\times \mathbb{R}_+.
	\end{aligned}
\end{equation}
%In order to absorb the quantity $-{\bf i}\sigma\omega^b(k)\hat{a}(k,\sigma,t)$ on the right hand side of \eqref{StartPoint0}, we set
%\begin{equation}
%	\label{HatA}
%	\alpha(k,\sigma,t) \ = \ \hat{a}(k,\sigma,t)e^{{\bf i}\sigma \omega^b(k)t}.
%\end{equation}
For  sake of simplicity, we also denote $\hat{a}(k,\sigma,t)$  as $\hat{a}_t(k,\sigma)$  
\begin{equation}
	\begin{aligned}\label{StartPoint}
			\partial _t \hat{a}_t(k,\sigma)\ & =  \ -{\bf i}\sigma\bar\omega(k)\hat{a}_t(k,\sigma) \mathrm{d}t\ - \ {\bf i}\sigma\lambda\sum_{\sigma_1,\sigma_2\in\{\pm1\}}\sum_{k_1,k_2\in\Lambda^*}\delta( \sigma k-\sigma_1k_1-\sigma_2k_2)\\
		& \ \   \times \mathcal{M}(k,k_1,k_2)\hat{a}_t(k_1,\sigma_1) \hat{a}_t(k_2,\sigma_2),\\
		\\
		\hat{a}_0(k,1) \ & = \ a_0(k), \ \ \ \ \forall (k,t)\in\Lambda^*\times \mathbb{R}_+.
	\end{aligned}
\end{equation}

By setting 
 $$f_{\lambda,D}(k,t)=\langle\hat{a}_t(k,-1),\hat{a}_t(h^{-1}k,1)\rangle,$$ and scaling $k\to hk$, in the kinetic limit of  $D\to \infty$, $\lambda\to 0$ and $t=\lambda^{-2}\tau=\mathcal{O}(\lambda^{-2})$, under suitable randomization of the system, we obtain \cite{staffilani2021wave}
$$\lim_{\lambda\to 0, D\to\infty} f_{\lambda,D}(k,\lambda^{-2}\tau) = f(k,\tau)$$
which solves the wave turbulence model

\begin{equation}\label{WaveTurbulence}
	\begin{aligned}
	&	\partial_\tau f(k,t) \  =  \ \mathcal C[f](k),\ \ \ \ f(k,0) \ = \ f_0(k),\ \  \forall k\in\mathbb{T}^d,\\
	&	\mathcal C[f](k)\  = \ \int_{\mathbb{T}^6}K(\omega,\omega_1,\omega_2)\delta(k-k_1-k_2)\delta(\omega-\omega_1-\omega_2)[f_1f_2-ff_1-ff_2]\mathrm{d}k_1\mathrm{d}k_2\\
		& - \ 2\int_{\mathbb{T}^6}K(\omega,\omega_1,\omega_2)\delta(k_1-k-k_2)\delta(\omega_1-\omega-\omega_2)[f_2f-ff_1-f_1f_2]\mathrm{d}k_1\mathrm{d}k_2,
	\end{aligned}
\end{equation} 
where $f=f(k),$ $f_1=f(k_1),$ $f_2=f(k_2),$ $\omega=\omega(k),$ $\omega_1=\omega(k_1),$ $\omega_2=f(k_2)$ and $K(k,k_1,k_2) \ = \ [8\omega(k)\omega(k_1)\omega(k_2)]^{-1}$.

\section{Collective behavior kinetic models of discrete non-local wave equations with long-range  interactions}\label{Collective}

The model \eqref{HamiltonODEs} is indeed a  chain of anharmonic
oscillators, in which the Hamiltonian H given by \eqref{Hamiltonian}. The collective behavior of chains of anharmonic
oscillators is a subject of growing interests in the physical community. For instance, in the case of the Fermi-Pasta-Ulam chains, the collective  behavior can be obtained via  a long-range interaction generalisation, in which the interactions are chosen to be non-local  (see for instance \cite{christodoulidi2016dynamics,christodoulidi2014fermi}). Inspired by this idea, we replace the local operator $\mathcal{L}$ by a non-local one. As an illustration, we consider $\mathcal{L}=(-\Delta)^{\alpha}$ with $0<\alpha<1$ and obtain the following discrete wave equation with long-range lattice interactions (see \cite{kirkpatrick2013continuum} for the same setting for the nonlinear Sch\"odinger equation)
\begin{equation}
	\label{CollectiveBehavior:1}\begin{aligned}
		\partial_{tt}\psi(x,t) \ = &  \ \ h^d\sum_{y\in  \Lambda, y\ne x} \frac{\psi(y,t)-\psi(x,t)}{|y-x|^{d+2\alpha}} \ - \ \lambda ({\psi(x,t)})^2, \\
		\psi(x,0) \ & = \ \psi_0(x), \ \partial_{t}\psi(x,0) \  = \ \psi_1(x), \ \forall (x,t)\in\Lambda\times \mathbb{R}_+.\end{aligned}
\end{equation}
Note that $h^d\sum_{y\in  \Lambda, y\ne x} \frac{\psi(y,t)-\psi(x,t)}{|y-x|^{d+2\alpha}} $ is the discretized version of the fractional Laplacian
\begin{equation}
	\label{CollectiveBehavior:1:a}\begin{aligned}
		(-\Delta \psi)^\alpha \ = &  -\ C_{d,\alpha}\int_{\mathbb{T}^d}\mathrm{d}y \frac{\psi(y,t)-\psi(x,t)}{|y-x|^{d+2\alpha}},\end{aligned}
\end{equation}
with $C_{d,\alpha}=\frac{4^\alpha \Gamma(d/2+\alpha)}{\pi^{\frac{d}{2}}|\Gamma(-\alpha)|}$.

In the scope of our paper, we restrict our considerations to the linear case $\lambda=0$
\begin{equation}
	\label{CollectiveBehavior:2}\begin{aligned}
		\partial_{tt}\psi(x,t) \ & =   \ \  h^d\sum_{y\in  \Lambda, y\ne x} \frac{\psi(y,t)-\psi(x,t)}{|y-x|^{d+2\alpha}}, \\
		\psi(x,0) \ & = \ \psi_0(x), \ \partial_{t}\psi(x,0) \  = \ \psi_1(x), \ \forall (x,t)\in\Lambda\times \mathbb{R}_+,\end{aligned}
\end{equation}
or equivalently, if we set $r_x(t)=\psi(x,t)$ and $v_x(t)=\partial_{t}\psi(x,t)$, the following system can be obtained (cf. \eqref{LatticeDynamicsSystem})
\begin{equation}
	\label{CollectiveBehavior:3}\begin{aligned}
		\partial_t{r}_x \ = & \ {v}_x,\ \ \ \ \ 
		\partial_t{v}_x \ =   \ \  h^d\sum_{y\in  \Lambda, y\ne x} \frac{r_y-r_x}{|y-x|^{d+2\alpha}}, \\
		r_x(0) \  = & \ \psi_0(x), \ v_x(0) \  = \ \psi_1(x), \ \forall (x,t)\in\Lambda\times \mathbb{R}_+.\end{aligned}
\end{equation}
Under suitable randomization of the initial conditions $\psi_0(x)$ and $\psi_1(x)$,
{this system of equations describes the long-range interactions of the lattice points $\Lambda$ and $$(\psi(x,t),\partial_{t}{\psi}(x,t))$$ represents the phase space position of the $x$-particle at time $t$.} Due to the long-range interactions between the $N=(2D+1)^d$ particles, a collective behavior dynamics is expected for \eqref{CollectiveBehavior:2}, similar to what happens for the Fermi-Pasta-Ulam chains with long-range interactions  \cite{christodoulidi2016dynamics,christodoulidi2014fermi}. 

Next, we will discuss ``mean-field limit'' of the above system by taking the limit $h\to 0$ or, equivalently $h^{-d}=N=(2D+1)^d\to \infty$. 
Formally, it is not a difficult task to derive the mean-field limit equation for Hamiltonian dynamics, see for instance \cite{ha2008particle} for the case of the Cucker-Smale model.  Rigorously, such derivations are challenging, especially when the interaction potentials are singular which is the current case.  We assume that the initial data are chosen in a way
that the empirical measure $N^{-1}\sum_{x\in\Lambda}\delta_{Q_x}\delta_{P_x}$ weakly converges in the limit $N\to\infty$ to the to the absolutely continuous
measure $g_0(r,v)\mathrm{d}r\mathrm{d}v$ with some smooth density $g_0(r,v)$. Here,  $r$ and $v$ are numbers in $\mathbb{R}$. We ask whether at some positive time $t>0$
the empirical measure $N^{-1}\sum_{x\in\Lambda}\delta_{Q_x(t)}\delta_{P_x(t)}$ weakly converges
to $g(r,v,t)\mathrm{d}r\mathrm{d}v$ with a density $g(r,v,t)$
satisfying some limiting evolution equation. Physically, the equation follows from the Liouville theorem, 
%which says that  a continuous medium in which each point moves under the action of an acceleration field behaves as an incompressible fluid. In our case, the form of our long-range interactions $h\sum_{y\in  \Lambda} \frac{Q_y-Q_x}{|y-x|^{1+2\alpha}}$ does not allow one to assume that the particles are identically the same as in the C-S case. We need to modify this argument as follows. 
assuming that the number of particles is  large enough such that it becomes meaningful to observe the distribution function
\begin{equation}
	\label{BBGKY1}
g^N \ =  \ g^N(t,({x,}r_x,v_x)_{x\in\Lambda}). 
\end{equation}
Defining the one-particle marginal distribution
\begin{equation}
	\label{BBGKY2}
	\rho^N\ (t,{x,}r_x,v_x)  \ =  \ \int_{\mathbb{R}^{(d+2)(|\Lambda|-1)}}\prod_{z\in\Lambda\backslash\{x\}}{\rm d}z\,\mathrm{d}r_z\,\mathrm{d}v_z\, g^N(t,({z},r_z,v_z)_{z\in\Lambda}), 
\end{equation}
where $|\Lambda|$ denotes the number of grid points. We now follow the BBGKY hierarchy to derive formally the kinetic description.  To this end denote $\Delta^h$ the discrete Laplacian in $d$-dimensions with $d\geq2$.  Set $\Phi^{h}_{\beta}(x) = (\Delta^{h})^{-1} |x|^{\beta}$ for $\beta\neq-2$ and $x\in\Lambda$.  Since $\Delta |\cdot|^{\beta+2} = (\beta+2)(d+\beta)|\cdot|^{\beta}$ we have that
\begin{equation*}
\Phi^{h}_{\beta}(x) \rightarrow c_{\beta,d}\,|x|^{\beta+2}\qquad\text{as}\qquad h\rightarrow0\,,\qquad c_{\beta,d}=\frac{1}{(\beta+2)(d+\beta)}\,,
\end{equation*}
in which the boundary condition of the problem $\Delta \Phi_{\beta} = |x|^{\beta}$ is chosen appropriately such that the discretised sequence $\{\Phi^{h}_{\beta}(x)\}_{h\in\Lambda}$ has the desired limit.  The Liouville equation reads, setting in the sequel $\beta=-d-2\alpha$,
\begin{equation}
	\label{BBGKY3}
\partial_tg^N \ + \ \sum_{x\in\Lambda}v_x\partial_{r_x}g^N  \ + \  h^d\sum_{x\in\Lambda}\partial_{v_x} \Big(\sum_{y\in  \Lambda, y\ne x} (r_y-r_x)(\Delta^{h}\Phi^{h}_{\beta})(y-x) g^N\Big) \ = \ 0\,.
\end{equation}
%which is equivalent to
%\begin{equation}
%	\label{BBGKY4}
%	\partial_tg^N \ + \ \sum_{x\in\Lambda}v_x\partial_{r_x}g^N  \ + \  \frac{1}{N^d}\sum_{x\in\Lambda}\partial_{v_x} \Big(\sum_{y\in  \Lambda, y\ne x}\frac{r_y-r_x}{|y-x|^{d+2\alpha}} g^N\Big) \ = \ 0.
%\end{equation}
We will now integrate both sides of \eqref{BBGKY3} with respect to $\mathrm{d}r_y\mathrm{d}v_y$, with $y\in\Lambda\backslash\{x\}$, to study the marginal distribution $\rho^N(x,r_x,v_x)$. Under the assumption that $g^N$ is rapidly decaying at infinity, the transport term in \eqref{BBGKY3} amounts to
\begin{multline}
	\label{BBGKY5}
	 \int_{\mathbb{R}^{(d+2)(|\Lambda|-1)}}\prod_{y\in\Lambda\backslash\{x\}}{{\rm d}y}\,\mathrm{d}r_y\mathrm{d}v_y \Big(\sum_{z\in\Lambda}v_z\partial_{r_z} g^N(t,(z,r_z,v_z)_{z\in\Lambda})\Big) \\ = \ v_x\partial_{r_x}\rho^N(t,x,r_x,v_x). 
\end{multline}
We next study the forcing term, which, by integration by parts reads
\begin{equation}
	\label{BBGKY6}\begin{aligned}
&	h^d\int_{\mathbb{R}^{(d+2)(|\Lambda|-1)}}\prod_{y\in\Lambda\backslash\{x\}}{{\rm d}y}\,\mathrm{d}r_y\mathrm{d}v_y\Big[ \sum_{z\in\Lambda}\partial_{v_z} \Big(\sum_{s\in  \Lambda, s\ne z}(r_s-r_z)(\Delta^{h}\Phi^{h}_{\beta})(s-z)g^N\Big)\Big] \\
= \ &	h^d\int_{\mathbb{R}^{(d+2)(|\Lambda|-1)}}\prod_{y\in\Lambda\backslash\{x\}}{{\rm d}y}\,\mathrm{d}r_y\mathrm{d}v_y \Big[\partial_{v_x} \Big(\sum_{s\in  \Lambda, s\ne x}(r_s-r_x)(\Delta^{h}\Phi^{h}_{\beta})(s-x)\,g^N\Big)\Big]. 
	 \end{aligned}
\end{equation}
We  move $\partial_{v_x}$ and the quantity $\sum_{s\in  \Lambda, s\ne x}\Phi^h_{\beta}(s-x)$ outside of the integral and obtain
\begin{equation}
	\label{BBGKY7}\begin{aligned}
		&	h^d\int_{\mathbb{R}^{(d+2)(|\Lambda|-1)}}\prod_{y\in\Lambda\backslash\{x\}}{{\rm d}y}\,\mathrm{d}r_y\mathrm{d}v_y \Big[\sum_{z\in\Lambda}\partial_{v_z} \Big(\sum_{y\in  \Lambda, y\ne z}(r_y-r_z) (\Delta^{h}\Phi^{h}_{\beta})(y-z)g^N\Big)\Big]\\
		= \ &	\partial_{v_x} \Big( h^d \sum_{s\in  \Lambda, s\ne x}(\Delta^{h}\Phi^h_{\beta})(s-x)\int_{\mathbb{R}^{(d+2)(|\Lambda|-1)}}\prod_{y\in\Lambda\backslash\{x\}}{{\rm d}y}\,\mathrm{d}r_y\mathrm{d}v_y (r_s-r_x)g^N\Big). 
	\end{aligned}
\end{equation}
We define the two-particle marginal function 
\begin{equation}
	\label{BBGKY8}
	\varrho^N (t,x,r_x,v_x,{y},r_y,v_y)  \ =  \ \int_{\mathbb{R}^{(d+2)(|\Lambda|-1)}}\prod_{z\in\Lambda\backslash\{x,y\}}{{\rm d}z}\,\mathrm{d}r_z\mathrm{d}v_z\, g^N(t,({z,}r_z,v_z)_{z\in\Lambda}), 
\end{equation}
and find
\begin{equation}
	\label{BBGKY9}\begin{aligned}
		&	h^d\int_{\mathbb{R}^{(d+2)(|\Lambda|-1)}}\prod_{y\in\Lambda\backslash\{x\}}{{\rm d}y}\,\mathrm{d}r_y\mathrm{d}v_y \sum_{z\in\Lambda}\partial_{v_z} \Big(\sum_{y\in  \Lambda, y\ne z}(r_y-r_z)(\Delta^h\Phi^h_\beta)(y-z) g^N\Big) \\
		= \ &	\partial_{v_x} \Big( h^d \sum_{s\in  \Lambda, s\ne x}(\Delta^h\Phi^h_{\beta})(s-x)\int_{\mathbb{R}^{2+d}}\mathrm{d}s\,\mathrm{d}r_s\,\mathrm{d}v_s \,(r_s-r_x)\varrho^N\Big). 
	\end{aligned}
\end{equation}
This leads to the following equation for the one-particle marginal function $\rho^N$
\begin{equation}
	\label{BBGKY10}\begin{aligned}
&	\partial_t\rho^N\ + \ v_x\partial_{r_x}\rho^N \\
& \ + \  \partial_{v_x} \Big( h^d \sum_{s\in  \Lambda, s\ne x}(\Delta^h\Phi^h_{\beta})(s-x)\int_{\mathbb{R}^{2+d}}\mathrm{d}s\,\mathrm{d}r_s\,\mathrm{d}v_s\, (r_s-r_x)\varrho^N\Big) \ = \ 0.	\end{aligned}
\end{equation}
Recall that $\Delta^{h}$ is a self-adjoint operator acting on the lattice location variable $s\in\Lambda$, consequently
\begin{align*}
&h^d \sum_{s\in  \Lambda, s\ne x}(\Delta^h\Phi^h_{\beta})(s-x)\int_{\mathbb{R}^{2+d}}\mathrm{d}s\,\mathrm{d}r_s\,\mathrm{d}v_s\, (r_s-r_x)\varrho^N \\
= \ & h^d \sum_{s\in  \Lambda, s\ne x}(\Delta^h\Phi^h_{\beta})(s-x)\int_{\mathbb{R}^{2+d}}\mathrm{d}s\,\mathrm{d}\tilde r\,\mathrm{d}\tilde v \,(\tilde r - r_x)\,\varrho^N\\
= \ & \int_{\mathbb{R}^{2+d}}\mathrm{d}s\,\mathrm{d}\tilde r\,\mathrm{d}\tilde v\, (\tilde r - r_x) \,\langle \Phi^{h}_{\beta}(\cdot-x),(\Delta^{h}\varrho^N)(t,x,r,v,\cdot,\tilde r,\tilde v)\rangle
\end{align*}
Passing to the mean-field limit $N\to\infty$, we obtain the one- and two-particle density functions (dropping the sub-$x$ notation) 
\begin{equation}
	\label{BBGKY11}\begin{aligned}
		&	\lim_{N\to\infty}\rho^N(t,x,r_x,v_x) \ = \ g(t,x, r,v),\\
		&    \lim_{N\to\infty}\varrho^N(t,x,r_x,v_x,s,\tilde{r},\tilde{v}) \ = \ \tilde{g}(t,x,r,v,s,\tilde{r},\tilde{v}),	\end{aligned}
\end{equation}
and the formal limit (note that $\beta+2>-d$)
\begin{equation*}
\langle \Phi^{h}_{\beta}(x-\cdot),(\Delta^{h}\varrho^N)(t,x,\cdot,r,v,\tilde r,\tilde v)\rangle \rightarrow c_{\beta,d}\int_{\mathbb{R}^{d}}|s-x|^{{\beta+2}}(\Delta_{s}\tilde{g})(t,x,r,v,s,\tilde{r},\tilde{v}){\rm d}s\,,
\end{equation*}
which leads to the mean-field equation 
\begin{equation}
	\label{BBGKY12}\begin{aligned}
		&	\partial_tg(t,x,r,v) \ + \ v\,\partial_{r}g(t,x,r,v)\\
		& \ + \  c_{\beta,d}\,\partial_{v} \Big(\int_{\mathbb{R}^{2}}\mathrm{d}\tilde r\,\mathrm{d}\tilde v\,(\tilde r-r)\int_{\mathbb{T}^d} \frac{\mathrm{d}s}{|s-x|^{d+2\alpha-2}} (\Delta_s\tilde{g})(t,x,s,r,v,\tilde{r},\tilde{v})\Big) \ = \ 0.	\end{aligned}
\end{equation}
If, in addition, we make the molecular chaos assumption  
\begin{equation}
	\label{BBGKY13}
	\tilde{g}(t,x,s,r,v,\tilde{r},\tilde v) \ = \ g(t,x,r,v)g(t,s,\tilde{r},\tilde{v}),
\end{equation}
then,
\begin{align*}
&c_{\beta,d}\int_{\mathbb{R}^{2}}\mathrm{d}\tilde r\,\mathrm{d}\tilde v\,(\tilde r-r)\int_{\mathbb{T}^d} \frac{\mathrm{d}s}{|s-x|^{d+2\alpha-2}} (\Delta_s\tilde{g})(t,x,s,r,v,\tilde{r},\tilde{v})\\
&=c_{\beta,d} \,g(t,x,r,v)\int_{\mathbb{R}^{2}}\mathrm{d}\tilde r\,\mathrm{d}\tilde v\,(\tilde r-r)\int_{\mathbb{T}^d} \frac{\mathrm{d}s}{|s-x|^{d+2\alpha-2}} (\Delta_s g)(t,s,\tilde{r},\tilde{v})\\
&=C_{d,\alpha}^{-1}\,g(t,x,r,v)\int_{\mathbb{R}^{2}}\mathrm{d}\tilde r\,\mathrm{d}\tilde v\,( r-\tilde r)(-\Delta_x)^{\alpha}g(t,x,\tilde r,\tilde v)\\
& = :C_{d,\alpha}^{-1}\,\Sigma_g(t,x,r)\,g(t,x,r,v)\,,
\end{align*}
and deduce from \eqref{BBGKY12} the Vlasov-type equation
\begin{equation}\label{Vlasov1}
\partial_tg(t,x,r,v) \ + \ v\partial_{r}g(t,x,r,v) \ + \ \mathbb{C}_1[g,g](t,x,r,v) \ = \ 0,
\end{equation}
with 
\begin{align*}
\mathbb{C}_1[g,g] \ :&= C_{d,\alpha}^{-1}\, \Sigma_g(t,x,r)\,\partial_{v} g(t,x,r,v)\\
\Sigma_g(t,x,r) &= (-\Delta_x)^{\alpha}\int_{\mathbb{R}^{2}}\mathrm{d}\tilde r\,\mathrm{d}\tilde v\,( r-\tilde r)\,g(t,x,\tilde r,\tilde v).
\end{align*}
Here $(-\Delta_x)^{\alpha}$ is the fractional Laplacian in the torus.  When the wave equation is set in the whole space we obtain, by an identical formal argument, the fractional Laplacian in $\mathbb{R}^d$.

%\textcolor{red}{do we need to rewrite the model?}
%\begin{equation}
%	\label{Vlasov3}\begin{aligned}
%		&	\partial_tg(t,x,r,v) \ + \ v\partial_{r}g(t,x,r,v) \ + \ \mathbb{C}_2[g](t,x,r,v) \ = \ 0,\\
%		&\mbox{ in which }\\
%		& \mathbb{C}_1[g,g] \ := -C_{d,\alpha}^{-1} \,\bigg(\int_{\mathbb{R}^{2}}\mathrm{d}\tilde r\mathrm{d}\tilde v(\tilde r-r)(-\Delta_x)^{\alpha}g(t,x,\tilde r,\tilde v)\bigg)\partial_{v}g(t,x,r,v)\,.\end{aligned}
%\end{equation}

%\bibliographystyle{plain}
%\bibliography{WaveTurbulence}

\def\cprime{$'$} \def\cprime{$'$} \def\cprime{$'$} \def\cprime{$'$}
\def\cprime{$'$} \def\cprime{$'$}

\end{document}